# Effects of Forecast Errors on Optimal Utilization in Aggregate Production Planning with Stochastic Customer Demand

**Klaus Altendorfer, Thomas Felberbauer, Herbert Jodlbauer**

Department of Operations Management, Upper Austria University of Applied Sciences

Wehrgrabengasse 1-3, A-4400 Steyr, Austria

klaus.altendorfer@fh-steyr.at

## Abstract

The hierarchical structure of production planning has the advantage of assigning different decision variables to their respective time horizons and therefore ensures their manageability. However, the restrictive structure of this top-down approach implying that upper level decisions are the constraints for lower level decisions also has its shortcomings. One problem that occurs is that deterministic mixed integer decision problems are often used for long term planning but the real production system faces a set of stochastic influences. Therefore, a planned utilization factor has to be included into this deterministic aggregate planning problem. In practise this decision is often based on past data and not consciously taken. In this paper the effect of long term forecast error on the optimal planned utilization factor is evaluated for a production system facing stochastic demand and the benefit of exploiting this decision's potential is discussed. Overall costs including capacity, backorder and inventory costs, are determined with simulation for different multi-stage and multi-item production system structures. The results show that the planned utilization factor used in the aggregate planning problem has a high influence on optimal costs. Additionally, the negative effect of forecast errors is evaluated and discussed in detail for different production system environments.

## 1    Introduction

In hierarchical production planning, mainly deterministic mixed integer decision problems are used for long term planning. In these deterministic problems, the objective is to minimize overall costs by fulfilling the customer demand and taking capacity constraints into consideration. This problem setup leads to utilizing nearly all available capacity, however, the real production system faces a set of stochastic influences which lead to poor service level and overall cost performance. Therefore, excess capacity has to be provided in the long term planning to react on these stochastic influences (Rajagopalan and Swaminathan 2001; Mieghem 1998). A practical approach to overcome this problem is to apply the average utilization from past time periods. However, this average utilization from past data is also linked to the past demand situation and therefore biased. This means that this approach ignores the conscious decision to treat the planned utilization factor as a decision variable which affects the overall costs. Another practical approach, observed in company projects conducted by the authors, is that the target utilization is used as strategic value. This strategic utilization value is employed without any awareness of the negative effects of a too high or low production system utilization.

A scientific approach to identify the optimal utilization factor applied in this study is to use simulation to quantify the benefits of consciously exploiting this decision variable. See Shafer and Smunt (2004) for a review of simulation applications and Tsai and Liu (2015) for an recent work in the area of operations management. Specifically, when the impact of forecast errors on costs is investigated extensive experimentation using simulation is a common method (Enns 2002; Jeunet 2006; Fildes and Kingsman 2010; Gansterer 2015). In this study, we investigate the effect of forecast error on the production system performance when the hierarchical production planning approach Manufacturing Resource Planning (MRP II) is applied. To support the general applicability of the simulation results and formulate some observations, a broad simulation study investigating different production system structures, demand patterns and cost rates is conducted.

Hierarchical production planning (HPP) models with decision levels differing in time horizon and decision variables are widely used (Hax and Meal 1975; Meal 1984; Schneeweiss 2003; Hopp and Spearman 2008), and its application has been proven to be very beneficial (Andersson, Axsäter, and Jönsson 1981; Gansterer 2015). The top-down approach of HPP causes that upper level decisions restrict the lower level decisions. Especially the MRP II concept, as presented in Hopp and Spearman (2008) or as discussed in Yeung, Wong, and Ma (1998), is a structured planning approach often implemented in enterprise resource planning systems and therefore frequently applied in practise. The MRP II planning approach hierarchy can be classified into three main planning levels, i.e. long-term planning, intermediate planning and short-term control. The long-term planning includes forecasting and aggregate production planning (APP) and determines appropriate levels of production, inventory and capacity (internal capacity, overtime and external capacity). Usually, APP models are solved with mixed integer linear programming (MILP) solvers (Hopp and Spearman 2008), and are implemented in the advance planning systems of commercial software vendors (Chern and Hsieh 2007; Kung and Chern 2009). The intermediate planning level computes the master production schedule (MPS) and the material requirements plan (MRP). For a detailed description of the MPS and the MRP algorithm we refer to Hopp and Spearman (2008) and Orlicky (1975). MRP is the most common material planning method, its position against other planning methods has strengthened since the 1990s (Jonsson and Mattsson 2006), and is still in focus of research (Enns 2001; Louly and Dolgui 2013). Finally, the lowest planning level, i.e. short-term control, performs the scheduling and dispatching tasks for production orders on the respective resources (Panwalkar and Iskander 1977; Framinan, Ruiz-Usano, and Leisten 2000; El-Bouri 2012).

There is a general agreement that uncertainty in forecasts is ubiquitous. In the literature different approaches are presented and investigated of how forecast error influences the production system. Ho and Ireland (1993) for example evaluate the impact of forecast error on the performance of a multi-product, multi-level MRP system by analysing planning methods nervousness. In a further study Ho and Ireland (1998) analyse the correlation between forecast errors and MRP nervousness in more detail. Using a simulation study they find that forecasting errors may not cause a higher degree of MRP nervousness and that the planning method parameter nervousness can be mitigated by using an appropriate lot sizing rule. Enns (2002) investigates the effects of forecast bias and demand uncertainty in a batch production system using an integrated MRP planning and execution test bed. Thereby, the use of inflated planned lead times and safety stock, to compensate for forecast error, is analysed. Another study of Jeunet (2006) analyses the impact of demand forecast errors on the cost performance of several lot-sizing techniques. The study shows that there is a nonlinear correlation of forecast errors with cost performance. Fildes and Kingsman (2010) investigate the effect of demand uncertainty and forecast error on costs and service levels for a MRP type manufacturing system. Ho and Ireland (2012) mitigate forecast errors by applying different lot sizing rules in enterprise resource planning controlled manufacturing systems. To evaluate the complex relationships within hierarchical production planning structures, a set of conference papers (Felberbauer, Altendorfer, and Hübl 2012; Felberbauer et al. 2013; Felberbauer and Altendorfer 2014) apply simulation to discuss different specific planning decisions. Felberbauer, Altendorfer, and Hübl (2012) investigate different machine allocation strategies when MRP and KANBAN controlled materials are produced on the same resources. Inspired by a simulation study for an automotive supplier, Felberbauer et al. (2013) compare the effect of different employee skill profiles, i.e. broader employee qualification, on inventory and service level. Based on the nervousness of MRP production orders, a customer information quality measure for hierarchical production planning is developed in Felberbauer and Altendorfer (2014). Using a similar planning setup, especially with respect to the APP, simulation is applied there to evaluate different customer order behaviours, i.e. "customer required lead time" and the "martingale model of forecast evolution", for one specific production system. Note that the forecast error and demand uncertainty effects, which are focussed in the current paper, are neglected in Felberbauer and Altendorfer (2014). In addition, the current study includes a comprehensive set of cost rates, different production system structures and demand patterns, as well as a set of practical implications.

The above reviewed literature shows that a lot of work has been done on the intermediate range planning methods MPS and MRP dealing with uncertainty by optimizing its parameters. See also Dolgui and Prodhon (2007) who present relevant literature on the parameterisation of MRP systems under demand and lead time uncertainties. The complex interdependence of the decisions in hierarchical production planning with uncertain information leads to a research gap concerning the influence of information uncertainties on APP decisions. Also, the surveys by Fleischmann and Meyr (2003), Kok and Fransoo (2003) and Missbauer and Uzsoy (2011) report the importance of integrating stochastic information into the decision process.

A recent study of Gansterer (2015) is a first step to close this research gap. Nevertheless, the main focus of Gansterer (2015) lies on the investigation of different forecasting techniques and neglects the discussion on how the deterministic APP model incorporates the stochastic nature of forecasts and demand which is addressed in the current study. This paper extends available literature by explicitly focusing on the effects of forecast error and demand uncertainty within a hierarchical production planning structure when applying a planned utilization factor in the APP.

The remainder of this article is structured as follows. In Section 2 the studied production systems and planning approach with the extension of the APP is explained. In section 3 the study design is introduced and detailed research questions are stated which are answered in the numerical study presented in Section 4. Some conclusions are provided in Section 5.

## 2    Studied production system and planning approach

As stated above, the APP problem is usually treated as deterministic (Schneeweiss 2003; Hax and Meal 1975) a production system faces a set of stochastic influences such as uncertainty in forecasts and demand. Therefore, in this study a planned utilization factor is integrated within the APP model and simulation is used to identify its optimal level to account for the forecast and demand uncertainty. In detail, the long term production plan is optimized (determining production amount and capacity levels) assuming perfect forecast information and is applied to a production system where the MRP II concept with MPS, MRP and dispatching functionalities is used. Forecast error, customer required lead time, order amount and order arrivals are modelled as random variables. The overall cost statement includes capacity and inventory costs, which are directly part of the APP problem, as well as backorder costs, which can only be evaluated by simulation of the stochastic system behaviour. Especially, the effect of forecast error, which incorporates the difference between deterministic monthly forecast and the stochastic realized monthly order quantities, on the optimal planned utilization factor, which is integrated in the APP model, is evaluated.

### 2.1    MRP II modeling

The standard MRP II planning approach consists of long, intermediate and short term planning levels which are calculated on rolling horizon basis. Figure 1 shows the planning model as applied in the current study. The forecast is assumed to be external information provided to the production system including monthly order amount per finished product. On the long term planning level, the APP (CPLEX© is used as optimization software) is conducted using the forecast, routing information, available shift models, planned utilization factor and current inventory levels as input. After optimization, the optimal production plan (planned monthly production quantities and planned inventory levels for each finished product) as well as the provided capacity, identified by the shift plan and external capacity, are the output of the MILP problem. The APP is conducted based on a rolling horizon planning approach 3 times a year.

On the operational level, single customer orders are generated for each finished product with a random order amount and a random due date. The order rate of these customer orders depends on the forecast. However, it is disturbed by a forecast error.

Both, customer orders and the optimal production plan according to the APP result are the input for the MPS where the monthly optimal production plan is disaggregated to a daily basis. The disaggregated optimal production plan is combined with the already available customer orders to generate the gross requirements for the MRP run. A cumulative calculation scheme is applied for the MPS calculation assuring that the maximum of cumulated customer demand and cumulated production program is satisfied with the cumulated MPS quantities. If there are planned inventory levels for a finished product from the APP solution, these planned inventory levels are also included in the MPS step. The gross requirements for finished products are the input for the MRP run with its steps: netting, lot sizing, backward scheduling and bill of material explosion. The output of the MRP run are the production orders for single materials with start date, end date and a production amount. The MRP and MPS functionality is conducted daily on a rolling horizon basis. The jobs are released continuously when their start date is reached and the necessary sub-materials are available. After job release, the production orders are produced according to their routing. Within the production process waiting production orders are sorted by the modified earliest due date (MEDD) dispatching rule according to their planned end date (Panwalkar and Iskander 1977). Due to the fact that in this study only the influence of stochastic demand and forecast error on the logistical performance is investigated all processing times are assumed to be deterministic and no setup time is applied. According to the suggestions of Hopp and Spearman (2008) for situations without setup times we apply the lot sizing policy lot-for-lot where the lot sizes are driven by specific customer orders. A detailed study of lot sizing effects is omitted. Machine availability is predefined for each month by the respective shift plan decision according to APP. The shop floor performance is measured as overall costs consisting of internal and external capacity costs, inventory costs and backorder costs.

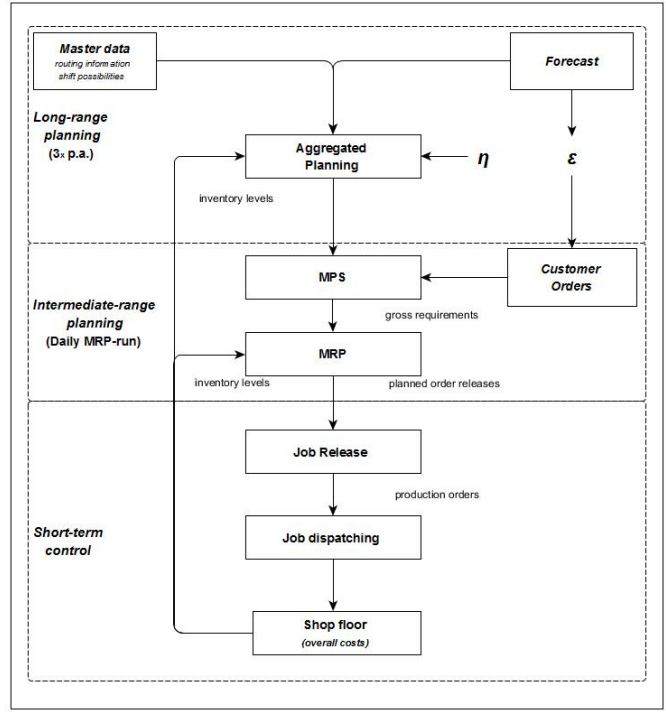

**Figure 1: Hierarchical production planning approach**

## 2.2 Details on APP problem studied

In this paper, a classical APP such as presented in Gansterer (2015), is extended by an additional decision on planned utilization (see also Felberbauer and Altendorfer 2014).

Table 1 shows all variables used in the MILP problem.

The APP is executed on a rolling horizon three times per year for a planning horizon of twelve months ($T$) which means that only the plan for the next four months is applied. Therefore the problem of zero inventory at the end of the planning horizon $T$ can be ignored.

The decision variables of the APP are: production program $x_{p,t}$ for finished product $p$ and month $t$, external capacity $e_{t,j}$ for resource $j$ in month $t$ and the shift model identified by the binary decision variable $w_{t,j,s} = 1$ if shift plan $s$, in time period $t$, and at machine $j$ is applied ($w_{t,j,s} = 0$ otherwise). Note that this external capacity, which ensures a feasible solution of the MILP problem, is assumed to be supplied by an external company which is able to provide the required technology of machine $j$. In the simulation model this external capacity is used in situations with high current system load. From the production program $x_{p,t}$, respective target inventory levels $l_{p,t}$ can be identified. These inventory levels state the planned pre-production of finished products, for example, to account for seasonality or for lower capacity in the next months. The following MILP problem is solved:

$$\sum_{t=1}^{T}\left(\sum_{j=1}^{J}\left(e_{t,j}c^e + \sum_{s=1}^{2}K_{t,j,s}c^i\,w_{t,j,s}\right) + \sum_{p=1}^{P}\tilde{c}^h\,l_{p,t}\right) \to \min_{\substack{x_{p,t}\\w_{t,j,s}\\e_{t,j}}} \tag{1}$$

subject to

$$\sum_{s=1}^{2}w_{t,j,s} = 1 \qquad\qquad \begin{aligned}&t = 1,\dots,T\\&j = 1,\dots,J\end{aligned} \tag{2}$$

$$\sum_{p=1}^{P}x_{p,t}a_{p,j} \le \eta\sum_{s=1}^{2}K_{t,j,s}w_{t,j,s} + e_{t,j} \qquad\qquad \begin{aligned}&t = 1,\dots,T\\&j = 1,\dots,J\end{aligned} \tag{3}$$

$$l_{p,t} = l_{p,t-1} + x_{p,t} - F_{p,t} \qquad\qquad \begin{aligned}&p = 1,\dots,P\\&t = 1,\dots,T\end{aligned} \tag{4}$$

$$l_{p,t}, x_{p,t} \ge 0 \qquad\qquad \begin{aligned}&p = 1,\dots,P\\&t = 1,\dots,T\end{aligned} \tag{5}$$

$$e_{t,j} \ge 0 \qquad\qquad \begin{aligned}&t = 1,\dots,T\\&j = 1,\dots,J\end{aligned} \tag{6}$$

$$w_{t,j,s} \in \{0,1\} \qquad\qquad \begin{aligned}&t = 1,\dots,T\\&j = 1,\dots,J\\&s = 1,\dots,2\end{aligned} \tag{7}$$

**Table 1: Symbols used (h=hour, d=day, m=month, CU= currency unit, pcs=pieces).**

| | | |
|---|---|---|
| $p = 1,\dots,P$ | Index of the finished products | 1 |
| $P$ | Index of the last finished product | 1 |
| $t = 1,\dots,T$ | Index of the time period month | m |
| $T$ | Index of the last time period of the planning horizon (12 month) | m |
| $\tau$ | Index of the time period day | d |
| $j = 1,\dots,J$ | Index of the machines | 1 |
| $J$ | Index of the last machine | 1 |
| $s \,\epsilon\, \{1,2\}$ | Index of the shift plan possibilities (10-shifts per week or 15-shifts per week) | 1 |
| $x_{p,t}$ | Decision variable production program of finished product $p$ in time period $t$ | pcs/m |
| $l_{p,t}$ | Inventory of finished product $p$ at the end of time period $t$ | pcs/m |
| $l_{p,0}$ | Initial inventory of finished product $p$ at the beginning of the aggregate planning | pcs/m |
| $w_{t,j,s}$ | Binary decision variable applying shift plan $s$ in time period $t$ at machine $j$ | 1 |
| $a_{p,j}$ | Required capacity for one piece of finished product $p$ on machine $j$ | h/pcs |
| $c^i$ | Internal capacity cost rate for one capacity unit at machine $j$ for shift plan $s$ | CU/h |
| $c^e$ | External capacity cost rate for one capacity unit at machine $j$ | CU/h |
| $\bar{c}^h$ | Inventory cost rate for one finished product $p$ per month t | CU/ m |
| $c^h$ | Inventory cost rate for one finished product $p$ per time period $\tau$ | CU/d |
| $c^b$ | Backorder cost rate for one finished product per time period $\tau$ | CU/d |
| $K_{t,j,s}$ | Available internal capacity of machine $j$ in time period $t$ applying shift plan $s$ | h/m |
| $e_{t,j}$ | Decision variable for the external capacity per machine $j$ in time period $t$ | h/m |
| $\eta$ | Planned utilization factor which defines the planned internal production system utilization | 1 |
| $\varepsilon_{p,t}$ | Forecast error which is identically independent normally distributed | pcs/m |
| $D_{p,t}$ | Real demand for finished product $p$ in time period $t$ | pcs/m |
| $o_p$ | The actual amount per order and finished product $p$ which is identically independent log-normally distributed | pcs |
| $\lambda_{p,t}$ | Arrival rate of finished product $p$ in time period $t$ | 1 |
| $\rho$ | Average annual capacity demand | shifts |
| $L$ | Customer required lead time which is identically independent log-normally distributed | d |
| $F_p$ | Annual forecast of finished product $p$ | pcs |
| $F_{p,t}$ | Forecast of finished product $p$ in time period $t$ | pcs/m |
| $\alpha$ | Forecast error parameter for the calculation of the forecast error deviation | 1 |

The objective function (1) minimizes the costs of internal and external capacity (at cost rate $c^i$ and $c^e$ for machine $j$ respectively) and the inventory costs (inventory cost rate per day $c^h$) which accrue by fulfilling the required demand. The main trade-off treated in the APP problem is between pre-production on stock and application of additional internal or external capacity assuming that the internal capacity cost rate is lower than the external one ($c^i \ll c^e$). Due to constraint (2) only one shift plan can be applied for one machine per time period $t$. Two possibilities for the shift plan decisions are available, which are the 10-shifts per week and the 15-shifts per week alternative. Constraints (3) ensure that the capacity needed to produce the optimal production program on the machines is lower than or equal to the available capacity of internal and external resources. The required capacity with respect to time period $t$ and machine $j$ is calculated using the processing time $a_{p,j}$ per finished product $p$ and machine $j$ and the respective optimal production program $x_{p,t}$. The external capacity is unlimited and the internal available capacity per machine $j$ is dependent on the applied shift plan and the planned utilization factor $\eta$. This factor is a real number between zero and one and defines the planned internal capacity utilization. Constraints (4) ensure that for the planning horizon $T$ the forecast $F_{p,t}$ per period $t$ is fulfilled by the optimal production program whereby it is also possible to use inventory as well.

Finally, the constraints (5)-(7) define the decision variables. As the external capacity ensures that the deterministic MILP problem is always feasible, no backorder costs are included in the cost statement. In a real production system (and also in the simulation model applied in this study) where the forecast is not exactly met and additional short term disturbances occur (random required lead time and order amount of single customer orders), not all orders will be on time. Note that a lower $\eta$ leads to higher optimal costs in the MILP problem (as more capacity has to be provided) but on the short term decision level, this decreased $\eta$ leads to lower backorder costs since additional capacity is available to react on the short term stochastic influences. As this effect cannot directly be included in the MILP model, the solution approach in this paper is to enumerate a set of possible utilization factor values and evaluate their expected overall costs by simulation.

*2.3    Production system structures and simulation model*

A production system with 6 machines, M1 to M6 and scenario specific routing, number of finished products as well as bills of materials is simulated using AnyLogic®. The modelled production systems follow either a flow shop or a job shop structure. The production system structures are streamlined versions of different production facilities operating in the automotive sector. These structures have been observed in a set of company projects. Table 2 shows the two investigated bill of material structures (for 4 and 8 finished products, material 10 - 17) and the two routing possibilities (job shop and flow shop). The abbreviation, for the scenario with 4 finished products and 4 sub materials, is "low number of products" and the abbreviation for the scenario with 8 finished products and 8 sub materials is "many products". Raw materials (material 100, 120 and 130) are assumed to be always available. In addition to the shop structure and the number of products, a constant and a seasonal forecast are investigated, whereby the seasonality is modelled by a sinusoidal function. These two demand patterns are again chosen based on empirical data from a set of company projects. Furthermore, also Rossi, Kilic, and Tarim (2015) recently discussed constant and seasonal demand in an inventory model investigation. Nevertheless, there are many further demand pattern alternatives which limit the findings of this study.

**Table 2: Bill of material and routing data.**

flow shop many products

| material | machine | component | # pieces |
|---|---|---|---|
| 10 | M5 | 20 | 1 |
| 11 | M5 | 20 | 1 |
| 12 | M6 | 21 | 1 |
| 13 | M6 | 21 | 1 |
| 14 | M5 | 22 | 1 |
| 15 | M5 | 23 | 1 |
| 16 | M6 | 32 | 1 |
| 17 | M6 | 33 | 1 |
| 20 | M3 | 30 | 1 |
| 21 | M4 | 31 | 1 |
| 22 | M3 | 32 | 1 |
| 23 | M4 | 33 | 1 |
| 30 | M1 | 100 | 1 |
| 31 | M2 | 110 | 1 |
| 32 | M1 | 120 | 1 |
| 33 | M2 | 130 | 1 |

job shop

| material | machine | component | # pieces |
|---|---|---|---|
| 10 | M5 | 20 | 1 |
| 11 | M5 | 20 | 1 |
| 12 | M6 | 21 | 1 |
| 13 | M6 | 21 | 1 |
| 14 | M5 | 22 | 1 |
| 15 | M5 | 22 | 1 |
| 16 | M6 | 23 | 1 |
| 17 | M6 | 23 | 1 |
| 20 | M3 | 30 | 1 |
| 21 | M4 | 31 | 1 |
| 22 | M6 M4 | 32 | 1 |
| 23 | M1 M4 | 33 | 1 |
| 30 | M1 | 100 | 1 |
| 31 | M2 | 110 | 1 |
| 32 | M1 M3 | 120 | 1 |
| 33 | M2 M3 | 130 | 1 |

flow shop low number of products

| material | machine | component | # pieces |
|---|---|---|---|
| 10 | M5 | 20 | 1 |
| 11 | M5 | 20 | 1 |
| 12 | M6 | 21 | 1 |
| 13 | M6 | 21 | 1 |
| 20 | M3 | 30 | 1 |
| 21 | M4 | 31 | 1 |
| 30 | M1 | 100 | 1 |
| 31 | M2 | 110 | 1 |

To account for the unbiased forecast error the random monthly customer demand $D_{p,t}$ is modelled as being the deterministic forecast value of this month $F_{p,t}$ disturbed by a random error term $\varepsilon_{p,t}$. This leads to $D_{p,t} = F_{p,t} + \varepsilon_{p,t}$. This forecast error $\varepsilon_{p,t}$ is an identically independent truncated-normal-distributed random variable with an expected value $\mathrm{E}[\varepsilon_{p,t}] = 0$ and variance $Var[\varepsilon_{p,t}] = (\alpha F_{p,t})^2$. The finished products $p$ are independently being forecasted. Therefore, the forecast error parameter $\alpha$, which is independent of time period $t$ and item $p$, defines the quality of the forecast. A limitation of this modelling assumption is that an unbiased forecast is assumed for which the true long term behaviour is known. The short term disturbances are included by a lognormal distributed amount per order $O_p$ (for finished product $p$) and a lognormal distributed customer required lead time $L$ (which has the same distribution for each finished product). Based on this order amount definition and the monthly demand, the order rate follows as: $\lambda_{p,t} = \frac{D_{p,t}}{E[O_p]} = \frac{F_{p,t} + \varepsilon_{p,t}}{E[O_p]}$. Note that in the simulation study the order rate $\lambda_{p,t}$ is adjusted monthly to account for the forecast error and the seasonality. In literature the above described order behaviour is known as "customer required lead time" order behaviour. The differences to another common order behaviour, i.e. the "martingale model of forecast evolution", is discussed Felberbauer and Altendorfer (2014).

The following example should clarify the relationship between the forecast, order rate and a customer order. Assuming a forecast value of $F_{p,t}$=1000 pcs per month, a forecast error value of $\varepsilon_{p,t} = 200$ and an expected value of a single customer order $\mathrm{E}[O] = 10$ would lead to an order rate of $\lambda_{p,t} = 120$. This means that the customer is going to order 120 times within this specific month where he requires for each of this orders an expected order amount of $\mathrm{E}[O] = 10$ and a random customer required lead time $L$.

The hierarchical production planning approach as described above is modelled with a simulation generator (Hübl et al. 2011; Felberbauer, Altendorfer, and Hübl 2012; Felberbauer and Altendorfer 2014), whereby the calculation of the APP is a function call within the simulation model. During simulation there is a cross data exchange between the simulation model and the MILP-Solver as described in previous section. The building blocks of the discrete event simulation model (built in AnyLogic®) are described in Hübl et al. (2011).

## 3    Study design

The following research questions are addressed in this study:

1) What is the impact of the planned utilization factor in a hierarchical production planning system on overall costs when forecasts are accurate or inaccurate and how do different demand or production system structures influence it? (Section 4.1)

2) What cost penalty results, if the planned utilization is not a decision variable but predefined at certain (past) values? (Section 4.2)

3) What influence does the level of inaccuracy of forecasts have on the optimal costs and optimal utilization? (Section 4.3)

4) How sensitive is this influence to: shop structure, number of products, demand pattern, capacity costs and backorder costs? (Section 4.4)

Section 4 presents the simulation study results, with detailed discussions of these research questions in the sub-sections. The following tables, Table 3 and Table 4, provide an overview over all parameters and scenarios applied to discuss the formulated questions.

**Table 3: Scenario specific parameters.**

Average shop load $\rho$ in shifts per day:    low … $\rho$=2.2 shifts/d

med … $\rho$=2.5 shifts/d

high … $\rho$=2.8 shifts/d

Internal $c^i$ and external $c^e$ capacity cost rates:    low … $c^i = 50$; $c^e = 100\ CU/h$

med … $c^i = 100$; $c^e = 200\ CU/h$

high … $c^i = 200$; $c^e = 400\ CU/h$

Backorder cost rate $c^b$ (independent of product)[1]:    low … $c^b = 9\ CU/d$

med … $c^b = 19\ CU/d$

high … $c^b = 99\ CU/d$

Forecast $F_{p,t}$ for constant demand $p\ \epsilon\ \{10,12,14,16\}$:    $F_{p,t} = 1000\ pcs/m$

Forecast $F_{p,t}$ for constant demand $p\ \epsilon\ \{11,13,15,17\}$:    $F_{p,t} = 1500\ pcs/m$

Forecast $F_{p,t}$ for seasonal demand $p\ \epsilon\ \{10,12,14,16\}$:    $F_{p,t} = 1000 + 500sin\left(2\pi\dfrac{t-5}{12}\right)$

$pcs/m$

Forecast $F_{p,t}$ for seasonal demand $p\ \epsilon\ \{11,13,15,17\}$:    $F_{p,t} = 1500 + 750sin\left(2\pi\dfrac{t-5}{12}\right)$

$pcs/m$

[1] Axsäter (2000) shows for streamlined inventory models that the low ratio leads to a target service level of 90%, the medium ratio leads to a target service level of 95% and the high ratio leads to a target service level of 99% which are feasible objectives for real production systems.

**Table 4: Parameters independent of scenario.**

Order amount finished product 10, 12, 14 and 16:    $E[O_p] = 10$ pcs; $Var[O_p] = 2.25$ pcs

Order amount finished product 11, 13, 15 and 17:    $E[O_p] = 15$ pcs; $Var[O_p] = 6.25$pcs

Customer required lead time for finished products:    $E[L] = 3$ d; $Var[L] = 3\ d$

Machine capacity for the 10-shifts per week plan:    $K_{t,j,1} = 320\ h$

Machine capacity for 15-shifts per week plan:    $K_{t,j,2} = 480\ h$

APP planning frequency:    3 times a year

APP Planning horizon:    $T = 12\ m$

MPS planning frequency:    daily

MPS planning horizon:    $60\ d$

MRP planning frequency:    daily

MRP planning horizon:    $30\ d$

MRP lotsizing policy for all materials:    lot-for-lot[1]

MRP planned lead time for all materials:    $2\ d$

MRP safety stock:    $ss_p = \dfrac{0.1F_p}{12}$

Materials inventory holding costs: $c_p^h = 0.5\ CU/d\quad \forall\ p\ \epsilon\ \{20,21,\ldots,33\}$

Finished products inventory holding costs: $c_p^h = 1\ CU/d\quad \forall\ p\ \epsilon\ \{10,11,\ldots,17\}$

Forecast error parameter: $\alpha = 0.05\kappa\ \forall\ \kappa\ \epsilon\ \{0,1,\ldots,10\}$[2]

Planned utilization values enumerated: $\eta = 0.5 + 0.02Q\ \forall\ Q\ \epsilon\ \{0,1,\ldots,25\}$

[1] see Hopp and Spearman (2008) for situations without setup times

[2] see Fildes and Kingsman (2010) for similar levels of forecast error uncertainty

The deterministic processing times are adjusted in each scenario to reach the predefined shop loads $\rho$ per machine $j$. For the constant demand pattern, average demand capacity per month and machine $j$ is 352h / 400h / 448h when shop load $\rho$ is 2.2/2.5/2.8, respectively. Note that for all scenarios all machines have exactly the same capacity load and that the utilization predefined by the planned utilization factor is also reached in the simulation model. For the flow shop structure with shop load $\rho$=2.5 this leads to a deterministic processing time $a_{p,j} \approx 10$ minutes in the scenario with a low number of products. In the scenario with many products the deterministic processing time is $a_{p,j} \approx 5$ minutes. For the scenarios with seasonal demand pattern the demand requirement per machine $j$ changes for each month following the sinusoidal function. Assuming a shop load of $\rho$=2.5 this means that demand starts with 400h in month 1, decreases to its lowest value of 200h in month 4, and increases again later to its peak value of 600h in month 10. MRP planning parameters for all materials have been identified in a preliminary study, again applying enumeration schemes to find appropriate values, and are not changed within the simulation study. In the further course of the paper the expression optimal overall costs and optimal planned utilization factor is a synonym for the best value found using the above described enumeration of planned utilization factors $\eta$. The numerical study, where a full factorial design is used, investigates 162 different scenarios: 3 shop structures (flow shop with many products, flow shop with low number of products and job shop with many products), 2 demand patterns, 3 shop loads, 3 capacity cost rates, 3 backorder cost rates. Combined with the 11 forecast error parameters this leads to 1,782 test instances for which the planned utilization factor is enumerated (26 possible values are tested according to $\eta = 0.5 + 0.02Q\ \forall\ Q\ \epsilon\ \{0,1,\ldots,25\}$) and therefore 46,332 simulation iterations are evaluated. Each iteration is replicated 10 times to account for stochastic variance which leads to 463,320 simulation runs. Four whole years are simulated whereby the first year is the warm up time of the simulation model and therefore excluded from the analysis. Customer orders that cannot be fulfilled at the end of the simulation time lead to extra costs in the form of penalties for their delay.

## 4    Simulation study results

The results presented in this section are coded according to the respective production system structures combined with the demand pattern: f_m_c is **f**low shop with **m**any products and **c**onstant demand (which is also called basic scenario), j_m_s **j**ob shop with **m**any products and **s**easonal demand or f_l_s **f**low shop with **l**ow number of products and **s**easonal demand. Please mention that in section 4.1 - 4.3 the cost rates and shop load parameters are always the medium setting for the discussion of the basic scenario.

### 4.1    Impact of Planned Internal Utilization Factor η on optimal costs and service level

The influence of the planned utilization factor $\eta$ for the basic scenario, with medium levels setting of cost rates and shop load parameters, is studied here in detail.

In Figure 2 the case without forecast errors ($\alpha = 0$) and in Figure 3 the case with a moderate forecast error ($\alpha = 0.25$) is presented. Furthermore, the figures show the single cost shares and the service level reached with respect to the planned utilization factor η.

For the basic scenario (f_m_c), which is the production system with flow shop, many products, constant demand, and medium levels of shop load $\rho$=2.5, capacity cost rate $c^i = 100/c^e = 200$ and backorder cost rate $c^b = 19$, Figure 2a shows that a rather high utilization factor of 98% is optimal for the situation without forecast error. This exhibits that the short term stochastic influences of order amount $O$ and customer required lead time $L$ for the production system assuming perfect forecast information and no further sources of uncertainty are very well compensated by the MPS and MRP functionality. The MPS function smoothes the gross requirements which are input for the MRP run. Additionally to the MPS, also the MRP-parameters planned lead time and safety stock lead to a certain amount of pre-production which compensates for some of the material shortages resulting from the stochastic demand uncertainty. The further sources of uncertainty linked to a real production process (stochastic processing times, setup times and machine failures) which usually lead to lower optimal utilization values (see, Jodlbauer and Altendorfer 2010) are not the focus of this paper. Nevertheless, the simulation model and the APP optimization process including the planned utilization factor can also be applied to a situation with these kinds of uncertainty. The detailed cost structure from Figure 2b shows that until a very high planned utilization, this factor only affects the external costs (which clearly decrease with higher internal utilization). In situations with high utilization values also backorder costs increase substantially.

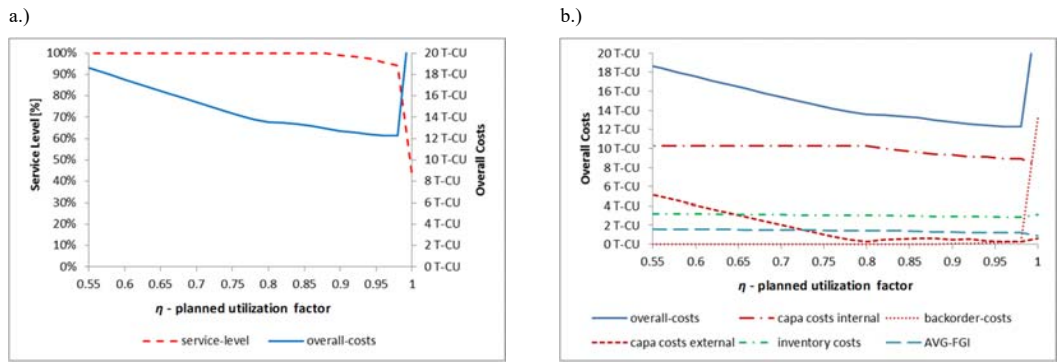

**Figure 2: Service level and costs with respect to the planned utilization factor $\eta$ for basic scenario without forecast error $\alpha=0$**

*A first **observation** from this setting, (cf research question 1) is that in the hierarchical production planning structure modelled (which is the MRP II structure) the MPS and MRP functionality mitigate the negative effects of short term customer demand uncertainties. However, also in this situation it is necessary to reduce the planned utilization factor below 100%.*

Figure 3 shows the results for a situation with inaccurate forecasts, $\alpha=0.25$, and already provides some further insights on the effect of such a forecast error. The optimal utilization, again for the basic scenario (f_m_c), is in this case 84 % and also the backorder costs become at utilization levels of above 85% already relevant. This means that for inaccurate forecasts a certain capacity has to be reserved to react on this inaccuracy. For both tested cases, the service level decreases very fast if the planned utilization factor is above the optimal value.

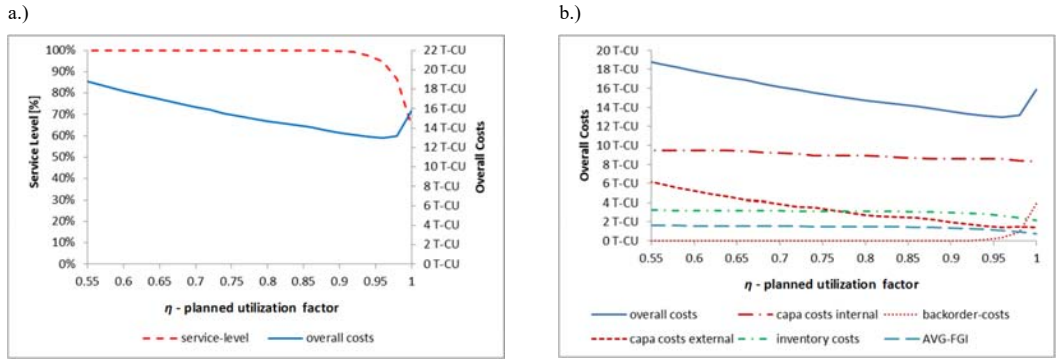

**Figure 3: Service level and costs with respect to the planned utilization factor η for basic scenario with forecast error α=0.25**

A further **observation** (cf. research question 1) from this setting is that in the hierarchical production planning structure modelled, it is necessary to include excess capacity in the APP optimization problem (for example with a planned utilization factor) to mitigate the effect that APP assumes a deterministic forecast but the realization of this forecast is stochastic.

Note that simulation is an appropriate tool to quantify the amount of excess capacity needed, however, the detailed results shown in Figure 3 are limited to the studied production system structure. To investigate the effects of different production systems and demand structures, the optimal utilization factors are evaluated by enumeration. Table 5 and Table 6 show the optimal planned utilization factors for the different scenarios, averaging over shop load $\rho$ and cost rates $c^i, c^e, c^b$, differentiating in cases without (see Table 5 where $\alpha = 0$) and cases with forecast error (see Table 6 where $\alpha = 0.25$), respectively. The "$\Delta$ in %" values show the difference to the basic scenario.

**Table 5: Optimal planned utilization η\* and costs without forecast errors α=0**

|  |  | scenario | | | | | |
|---|---|---|---|---|---|---|---|
|  |  | f  m  c\* | f  m  s | f  l  c | f  l  s | j  m  c | j  m  s |
| utilization | optimal | 96% | 96% | 92% | 93% | 96% | 96% |
|  | $\Delta$ in % | 0 | 0.00% | -4.17% | -3.12% | 0.00% | 0.00% |
| costs | optimal | 13,637 | 14,590 | 12,521 | 13,804 | 13,550 | 14,525 |
|  | $\Delta$ in % | 0 | 6.99% | -8.18% | 1.22% | -0.64% | 6.51% |

*f* ... flow shop          *j* ... job shop
*m* ... many products      *l* ... low number of products
*c* ... constant demand pattern    *s* ... seasonal demand pattern
\* cost and utilization Δ-values are compared to this scenario

**Table 6: Optimal planned utilization η\* and costs with forecast error α=0.25**

| | | scenario | | | | | |
|---|---|---|---|---|---|---|---|
| | | f_m_c* | f_m_s | f_l_c | f_l_s | j_m_c | j_m_s |
| utilization | optimal | 85% | 83% | 79% | 76% | 85% | 84% |
| | Δ in % | 0 | -2.35% | -7.06% | -10.59% | 0.00% | -1.18% |
| costs | optimal | 15,331 | 16,712 | 14,713 | 16,234 | 15,202 | 16,567 |
| | Δ in % | 0 | 9.01% | -4.03% | 5.89% | -0.84% | 8.06% |

\* cost and utilization Δ-values are compared to this scenario

The results of both Table 5 and Table 6, indicate that a lower number of products lead to lower costs and lower optimal planned utilization. The lower optimal planned utilization results from a balancing effect of forecast errors that more products have. The counterintuitive result of lower costs is linked to the cost structure in the simulation model. As backorder and holding costs are not dependent on the processing times the lower number of products, which leads to less throughput in pieces, leads to lower costs. In section 4.4, the effects of a lower number of products are further discussed.

*An **observation** (cf. research question 1) concerning seasonal demand pattern is that for both cases, without and with forecast error, the seasonal demand leads to higher costs in comparison to the constant demand.*

Intuitive reasons for this behaviour are higher external capacity costs in the peak months and additional inventory costs for pre-production. A further interesting result concerning seasonality is that for the situation without forecast error, the optimal planned utilization is nearly equal to the constant demand scenarios. When a forecast error occurs (see Table 6) the seasonal demand leads to a lower optimal planned utilization which may be linked to the higher absolute effect the forecast error has on the peak months.

The results from Table 5 and Table 6 show that there is nearly no difference between job shop (scenario **j**_m_c and **j**_m_s) and flow shop structure (scenario **f**_m_c and **f**_m_s). As in this study the difference between job shop and flow shop are only the routing paths, this result indicates that problems in job shop production systems are not mainly triggered by the more complex routing of products but have some other reasons. A discussion of these effects is left to further research.

*4.2    Cost penalty for ignoring planned utilization decision variable*

As motivated in the introduction, practise in hierarchical production planning is often to set the planned utilization at the average level from past periods, e.g. the last year, which ignores the cost reduction potential of treating $\eta$ as a decision variable. In this subsection the cost effect of this often unconsciously taken decision is discussed. The following Figure 4 shows the percentage in cost increase if the planned utilization factor is above or below the optimal value. In detail, for each of the 162 scenarios the optimal costs are compared with the costs reached when the planned utilization factor deviates from the optimum. The values shown in Figure 4 are the average, minimum and maximum percentage of cost increase. Figure 4a shows the case without forecast error $\alpha=0$ and Figure 4b the values for $\alpha=0.25$.

For the situation without forecast error, the results show that assuming a utilization value from past data which is too low, leads to a slight cost increase (only 1.35% average costs increase if utilization is set 2% below the optimal value) but higher utilization values lead to very high cost increases. This very fast increase in costs for high utilization values is based on the situation that without forecast error only little excess capacity is needed in the current setting to react on customer demand uncertainty (optimal planned utilization factor for basic scenario is 98%). Therefore the optimal planned utilization factor is near 100% but a planned utilization of exactly 100% is still far more expensive than a few percentage points below.

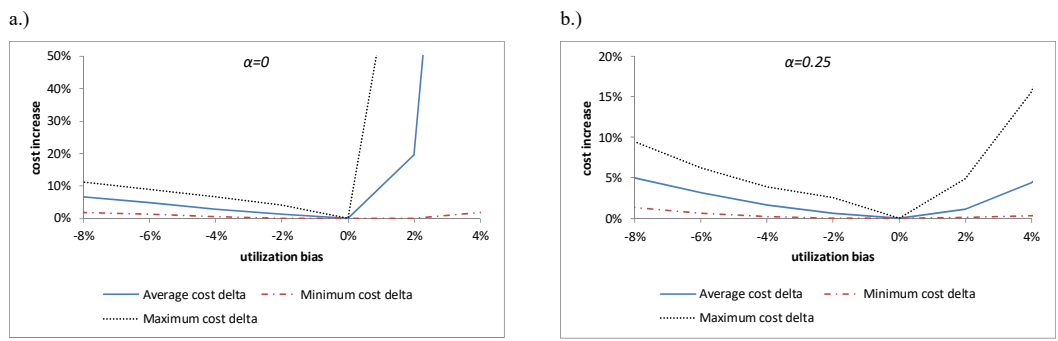

**Figure 4: Cost effect of wrong planned utilization factor $\eta$.**

Figure 4b shows that the relative cost increase for the situation with forecast error $\alpha=0.25$ is lower than in the situation without forecast error, which results from an average planned utilization of 85% to be optimal in this setting. Therefore, an increase by 2%, which leads to an average utilization of 87%, has less relative effect. Nevertheless, the absolute cost value reached, which is depicted in Table 5 and Table 6, is for the situation without forecast error 13,637 which is less than the 15,331 for the situation with forecast error.

Figure 4b shows that a too high planned utilization value leads also in this case to a considerable cost increase. On average a planned utilization value of 4% above the optimum leads to a cost increase of about 4%. For both situations the cost risk for too low planned utilization values is lower but still considerable. Looking at the maximum values of cost increases shows that for some single scenarios this decision variable has an even higher effect.

*This leads to the **observation** (cf. research question 2) that a high optimization potential is wasted if the planned utilization decision is not consciously taken. Furthermore, higher target utilizations (which could be a management tool) have a high risk of damaging the production system performance.* This second risk of too high target utilization values in production systems is based on higher backorder costs resulting from lower capacity provided.

### 4.3 Influence of forecast error α on optimal costs and planned utilization

Based on the simulation results for each test instance (where the optimal planned utilization factor is already identified), Figure 5 shows the average results over all test instances for the basic case f_m_c (**f**low shop, **m**any products and **c**onstant demand) and the f_l_c scenario (**f**low shop, **l**ow number of products and **c**onstant demand).

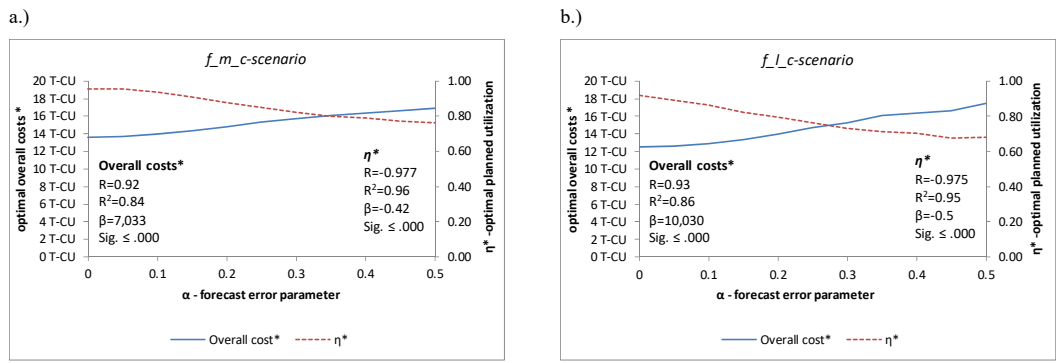

**Figure 5: Optimal costs and optimal planned utilization factor η\* with respect to the forecast error parameter α.**

In Figure 5, the optimal planned utilization factor is shown with respect to the forecast error and also the optimal overall costs are depicted. Furthermore, a linear regression analysis is conducted for each single test instance and the average values over all test instances are presented here. The finding already motivated by the results in section 4.1 is that an increase in forecast error leads to an increase of optimal overall costs. For flow shop production structure, many products and constant demand (f_m_c-scenario) an improvement of forecast quality from α=0.3 to α=0.1 leads to an decrease of overall cost by ≈11%, whereby the planned utilization could be increased from 82% to 94%. For the same production structure and demand pattern but a lower number of products (f_l_c-scenario) the cost reduction potential is ≈16% and the internal utilization level could be increased from 76% to 89%. Analysing the slopes $\beta$ and intercepts of the linear regression models for the f_m_c and f_l_c-scenario leads us to the finding that an improvement of forecast quality by α=-0.1 improves costs in average by 5% / 8%, respectively.

An **observation** (cf. research question 3) from these results is that higher long term information uncertainty, which leads to forecast errors, has to be compensated by a higher capacity provided by the aggregate planning which leads to significantly higher costs.

This finding is in line with Jeunet (2006) who also discusses the impact of demand forecast errors on the cost performance of several lot-sizing techniques in a multi-level context. Other studies show that there are non-linear dependencies of forecast errors and key performance indicators. Fildes and Kingsman (2010) found in a MRP controlled manufacturing system non-linear relationships of unit costs and demand uncertainty with respect to certain service level constraints. Nevertheless, Lee and Adam (1986) show that for a complex MRP controlled production system, which applies forecasts, the rationality of higher forecast error resulting in higher total cost may not be guaranteed. In results from Figure 5a and 5b, the introduction of the planned utilization factor and its adaption with respect to the unbiased forecast error lead to a nearly linear dependency of forecast error and costs.

The comparison between the optimal utilization in Figure 5a and Figure 5b further indicates that a lower number of products also means a higher capacity to be provided and shows that it holds for a broad range of forecast error situations.

An **observation** (cf. research question 3) is therefore, that hierarchical production planning systems facing unbiased (expected values 0) and uncorrelated (independent between items) forecast errors with more different products are more robust with respect to the forecast error influence.

This intuitive result is interesting as it shows that the hierarchical planning structure supports the underlying statistics (coefficient of variation of the sum of $n$ independent random variables reduces with respect to $n$) and that also this negative long term stochastic influence can be affected by the number of products. The higher average $\beta$-value in Figure 5b further shows that also the slope of the planned utilization decrease is for lower number of products steeper.

*4.4    Sensitivity analysis and broad numerical study*

To identify the sensitivity of the results with respect to shop structure, number of products, demand pattern, capacity costs and backorder costs, the following tables provide aggregated results of the full factorial simulation study design (only selected α-values are presented). The average results concerning optimal planned utilization factor and overall costs for the basic scenario f_m_c (flow shop, many products and constant demand pattern) over all respective test instances with the specified parameters are presented in Table 7. This shows that for a medium shop load level $\rho$=2.5 and a moderate forecast error parameter of α=0.2 the optimal planned utilization factor $\eta$* is 0.87 and the optimal overall costs are 14,805 for the f_m_c-scenario. Note that the optimal planned utilization factor $\eta$* of 0.87 is the average value over all capacity cost rate and backorder cost rate levels (low, med, and high). The column labelled with Ø reports the average value of $\eta$* and optimal overall costs over the forecast error parameter α with respect to the shop load $\rho$, capacity cost rates $c^i/c^e$ and the backorder cost rate $c^b$. The sub table rows indicated with Ø report the average value $\eta$* and optimal overall costs with respect to a specific forecast error parameter α.

**Table 7: Optimal planned utilization $\eta^*$ and overall costs– sensitivity analysis for basic scenario f_m_c.**

optimal planned utilization $\eta^*$

| | | a - forecast error parameter | | | | | | |
|---|---|---|---|---|---|---|---|---|
| | | 0 | 0.1 | 0.2 | 0.3 | 0.4 | 0.5 | Ø |
| $\rho$ | 2.2 | 0.95 | 0.93 | 0.87 | 0.81 | 0.78 | 0.76 | **0.85** |
| | 2.5 | 0.96 | 0.94 | 0.87 | 0.83 | 0.79 | 0.76 | **0.86** |
| | 2.8 | 0.96 | 0.94 | 0.89 | 0.82 | 0.80 | 0.77 | **0.86** |
| | **Ø** | **0.96** | **0.94** | **0.88** | **0.82** | **0.79** | **0.76** | **0.86** |

| | | a - forecast error parameter | | | | | | |
|---|---|---|---|---|---|---|---|---|
| | | 0 | 0.1 | 0.2 | 0.3 | 0.4 | 0.5 | Ø |
| $c^i/c^e$ | low | 0.96 | 0.93 | 0.87 | 0.80 | 0.77 | 0.74 | **0.85** |
| | med | 0.96 | 0.93 | 0.88 | 0.82 | 0.79 | 0.76 | **0.86** |
| | high | 0.96 | 0.94 | 0.89 | 0.84 | 0.81 | 0.78 | **0.87** |
| | **Ø** | **0.96** | **0.94** | **0.88** | **0.82** | **0.79** | **0.76** | **0.86** |

| | | a - forecast error parameter | | | | | | |
|---|---|---|---|---|---|---|---|---|
| | | 0 | 0.1 | 0.2 | 0.3 | 0.4 | 0.5 | Ø |
| $c^b$ | low | 0.98 | 0.96 | 0.91 | 0.86 | 0.83 | 0.81 | **0.89** |
| | med | 0.96 | 0.94 | 0.90 | 0.83 | 0.80 | 0.78 | **0.87** |
| | high | 0.93 | 0.91 | 0.83 | 0.78 | 0.74 | 0.70 | **0.81** |
| | **Ø** | **0.96** | **0.94** | **0.88** | **0.82** | **0.79** | **0.76** | **0.86** |

optimal overall costs*

| | | a - forecast error parameter | | | | | | |
|---|---|---|---|---|---|---|---|---|
| | | 0 | 0.1 | 0.2 | 0.3 | 0.4 | 0.5 | Ø |
| $\rho$ | 2.2 | 12,319 | 12,823 | 13,668 | 14,405 | 14,834 | 15,240 | **13,881** |
| | 2.5 | 13,732 | 14,163 | 14,805 | 15,394 | 16,168 | 16,650 | **15,152** |
| | 2.8 | 14,859 | 14,977 | 15,868 | 17,373 | 18,187 | 18,898 | **16,694** |
| | **Ø** | **13,637** | **13,987** | **14,780** | **15,724** | **16,396** | **16,929** | **15,242** |

| | | a - forecast error parameter | | | | | | |
|---|---|---|---|---|---|---|---|---|
| | | 0 | 0.1 | 0.2 | 0.3 | 0.4 | 0.5 | Ø |
| $c^i/c^e$ | low | 7,564 | 7,812 | 8,269 | 8,743 | 9,064 | 9,365 | **8,470** |
| | med | 12,162 | 12,470 | 13,216 | 14,066 | 14,667 | 15,149 | **13,622** |
| | high | 21,183 | 21,679 | 22,856 | 24,362 | 25,457 | 26,274 | **23,636** |
| | **Ø** | **13,637** | **13,987** | **14,780** | **15,724** | **16,396** | **16,929** | **15,242** |

| | | a - forecast error parameter | | | | | | |
|---|---|---|---|---|---|---|---|---|
| | | 0 | 0.1 | 0.2 | 0.3 | 0.4 | 0.5 | Ø |
| $c^b$ | low | 13,417 | 13,710 | 14,289 | 15,148 | 15,680 | 16,076 | **14,720** |
| | med | 13,538 | 13,908 | 14,651 | 15,542 | 16,167 | 16,688 | **15,082** |
| | high | 13,954 | 14,343 | 15,400 | 16,481 | 17,342 | 18,024 | **15,924** |
| | **Ø** | **13,637** | **13,987** | **14,780** | **15,724** | **16,396** | **16,929** | **15,242** |

Ø ... average value    $c^i/c^e$ ... capacity cost rates internal / external

$\rho$ ... shop load    $c^b$ ... backorder cost rate

In addition to the findings from section 4.3, the results in Table 7 show that there is only a slight influence of shop load $\rho$ on the optimal planned utilization but the overall costs increase significantly. Furthermore, a positive relation between the planned utilization factor and the increasing capacity cost rate $c^i/c^e$ can be observed for high $\alpha$ values whereby also overall costs increase. This finding supports previous literature stating that the more expensive capacity is, the higher the production system utilization should be (see, Jodlbauer and Altendorfer 2010). Based on this finding, an intuitive result is that a higher backorder cost rate $c^b$ leads to lower optimal planned utilization factor and higher overall costs.

The following tables, Table 8 and Table 9, show the aggregated results of a linear regression analysis which has been performed for each of the 162 scenarios to identify the relationship between forecast error $\alpha$ and the optimal planned utilization $\eta^*$ as well as between forecast error $\alpha$ and the optimal overall costs. The results are grouped by the different production system structures investigated. For the abbreviations of the scenario, we refer to the legend below Table 8. The presented $R$-values (correlation coefficient) and $\beta$-values (slope of the relationship) in the tables are the average values over the respective scenarios and parameter combinations. The "%" value describes the difference measured in percentage points between no forecast error ($\alpha=0$) and a high forecast error ($\alpha=0.5$). All correlations are significant at a level above 99.99% for each single scenario.

**Table 8: Aggregated results of the linear regression analysis for optimal planned utilization factor $\eta*$.**

Correlation coefficients for optimal utilization $\eta *$

| | | production structure and demand scenarios | | | | | | | | | | | | | | | | | |
|---|---|---|---|---|---|---|---|---|---|---|---|---|---|---|---|---|---|---|---|
| | | f m c | | | f m s | | | f l c | | | f l s | | | j m c | | | j m s | | |
| | | R | β | % | R | β | % | R | β | % | R | β | % | R | β | % | R | β | % |
| $\rho$ shop load | 2.2 | -0.97 | -0.43 | -21% | -0.97 | -0.45 | -23% | -0.97 | -0.50 | -26% | -0.97 | -0.60 | -32% | -0.97 | -0.41 | -20% | -0.97 | -0.42 | -22% |
| | 2.5 | -0.98 | -0.42 | -21% | -0.97 | -0.41 | -21% | -0.98 | -0.46 | -26% | -0.97 | -0.56 | -30% | -0.97 | -0.41 | -20% | -0.97 | -0.39 | -20% |
| | 2.8 | -0.98 | -0.42 | -20% | -0.97 | -0.43 | -21% | -0.98 | -0.55 | -29% | -0.96 | -0.53 | -28% | -0.98 | -0.39 | -19% | -0.97 | -0.40 | -20% |
| | Ø | **-0.98** | **-0.42** | **-20%** | **-0.97** | **-0.43** | **-22%** | **-0.98** | **-0.50** | **-27%** | **-0.96** | **-0.50** | **-30%** | **-0.98** | **-0.40** | **-19%** | **-0.97** | **-0.41** | **-21%** |

| | | production structure and demand scenarios | | | | | | | | | | | | | | | | | |
|---|---|---|---|---|---|---|---|---|---|---|---|---|---|---|---|---|---|---|---|
| | | f m c | | | f m s | | | f l c | | | f l s | | | j m c | | | j m s | | |
| | | R | β | % | R | β | % | R | β | % | R | β | % | R | β | % | R | β | % |
| $c^i/c^e$ capacity cost rates | low | -0.98 | -0.46 | -22% | -0.97 | -0.45 | -24% | -0.98 | -0.58 | -32% | -0.97 | -0.65 | -34% | -0.97 | -0.44 | -21% | -0.97 | -0.43 | -22% |
| | med | -0.98 | -0.42 | -20% | -0.98 | -0.44 | -22% | -0.97 | -0.49 | -26% | -0.97 | -0.60 | -31% | -0.98 | -0.41 | -20% | -0.98 | -0.41 | -21% |
| | high | -0.97 | -0.38 | -18% | -0.97 | -0.39 | -20% | -0.98 | -0.44 | -23% | -0.95 | -0.43 | -24% | -0.97 | -0.37 | -17% | -0.97 | -0.38 | -19% |
| | Ø | **-0.98** | **-0.42** | **-20%** | **-0.97** | **-0.43** | **-22%** | **-0.98** | **-0.50** | **-27%** | **-0.96** | **-0.50** | **-30%** | **-0.98** | **-0.40** | **-19%** | **-0.97** | **-0.41** | **-21%** |

| | | production structure and demand scenarios | | | | | | | | | | | | | | | | | |
|---|---|---|---|---|---|---|---|---|---|---|---|---|---|---|---|---|---|---|---|
| | | f m c | | | f m s | | | f l c | | | f l s | | | j m c | | | j m s | | |
| | | R | β | % | R | β | % | R | β | % | R | β | % | R | β | % | R | β | % |
| $c^b$ backorder cost rate | low | -0.98 | -0.37 | -17% | -0.97 | -0.38 | -19% | -0.98 | -0.43 | -22% | -0.96 | -0.50 | -27% | -0.97 | -0.35 | -16% | -0.97 | -0.36 | -18% |
| | med | -0.98 | -0.41 | -20% | -0.98 | -0.42 | -21% | -0.97 | -0.49 | -25% | -0.97 | -0.57 | -30% | -0.98 | -0.39 | -19% | -0.97 | -0.39 | -20% |
| | high | -0.98 | -0.48 | -25% | -0.97 | -0.49 | -26% | -0.98 | -0.59 | -33% | -0.96 | -0.61 | -33% | -0.98 | -0.47 | -24% | -0.97 | -0.41 | -24% |
| | Ø | **-0.98** | **-0.42** | **-20%** | **-0.97** | **-0.43** | **-22%** | **-0.98** | **-0.50** | **-27%** | **-0.96** | **-0.50** | **-30%** | **-0.98** | **-0.40** | **-19%** | **-0.97** | **-0.41** | **-21%** |

\* Sig. for $R$ is $\leq 0.000$

| | | |
|---|---|---|
| $f$ ... flow shop | $j$ ... job shop | $\emptyset$ ... average value | $c^i/c^e$ ... capacity cost rates internal / external |
| $m$ ... many products | $l$ ... low number of products | $\rho$ ... shop load | $c^b$ ... backorder cost rate |
| $c$ ... constant demand pattern | $s$ ... seasonal demand pattern | | |

An **observation** from Table 8 (cf. research question 4) is that the planned utilization factor decrease with increasing forecast error is found for any production system and cost structure. Furthermore, lower capacity cost rates and higher backorder cost rates lead for any production system and cost structure to a steeper planned utilization factor decrease.

The comparison between the production system structures shows, that especially the scenarios with a low number of products (labelled as f_l_c and f_l_s in Table 8) have a significantly higher cost increase when forecast error increases which supports the respective observation of section 4.3.

Therefore the general **observation** (cf. research question 4) related to forecast error variance is that companies facing high forecast error variances can achieve a significant cost reduction without traditional investments by improving customer collaboration and therefore getting more accurate forecast information, which may also imply some investment in information systems.

This finding is in line with traditional operations management literature (Wijngaard and Karaesmen 2007; Altendorfer and Minner 2011). Empirical evidence on the benefits through information sharing with substantial improvements in forecast accuracy is presented in Cui et al. (2015) and Trapero, Kourentzes, and Fildes (2012). In the current study this finding is extended to the accuracy of long term forecasts.

**Table 9: Aggregated results of the linear regression analysis for the optimal overall costs.**

**Correlation coefficients for overall costs***

| | | production structure and demand scenarios | | | | | | | | | | | | | | | | |
|---|---|---|---|---|---|---|---|---|---|---|---|---|---|---|---|---|---|---|
| | | f m c | | | f m s | | | f l c | | | f l s | | | j m c | | | j m s | | |
| ρ shop load | | R | β | % | R | β | % | R | β | % | R | β | % | R | β | % | R | β | % |
| | 2.2 | 0.93 | 5,948 | 24% | 0.88 | 6,533 | 26% | 0.93 | 7,277 | 33% | 0.91 | 8,161 | 35% | 0.92 | 5,569 | 22% | 0.87 | 6,131 | 23% |
| | 2.5 | 0.92 | 6,081 | 21% | 0.87 | 6,406 | 22% | 0.92 | 9,274 | 36% | 0.90 | 9,172 | 31% | 0.90 | 5,682 | 19% | 0.84 | 5,978 | 19% |
| | 2.8 | 0.91 | 9,071 | 27% | 0.87 | 7,558 | 23% | 0.94 | 13,539 | 48% | 0.89 | 10,182 | 31% | 0.90 | 8,580 | 25% | 0.86 | 7,228 | 22% |
| | **Ø** | **0.92** | **7,033** | **24%** | **0.87** | **6,832** | **23%** | **0.93** | **10,030** | **39%** | **0.70** | **13,368** | **35%** | **0.91** | **6,610** | **22%** | **0.86** | **6,446** | **21%** |

| | | production structure and demand scenarios | | | | | | | | | | | | | | | | |
|---|---|---|---|---|---|---|---|---|---|---|---|---|---|---|---|---|---|---|
| | | f m c | | | f m s | | | f l c | | | f l s | | | j m c | | | j m s | | |
| $c^i/c^e$ capacity cost rates | | R | β | % | R | β | % | R | β | % | R | β | % | R | β | % | R | β | % |
| | low | 0.91 | 3,711 | 24% | 0.85 | 3,517 | 22% | 0.93 | 5,396 | 44% | 0.91 | 5,124 | 39% | 0.90 | 3,512 | 22% | 0.84 | 3,416 | 21% |
| | med | 0.92 | 6,366 | 25% | 0.88 | 6,309 | 25% | 0.92 | 9,082 | 41% | 0.91 | 9,069 | 38% | 0.91 | 5,943 | 22% | 0.86 | 5,844 | 22% |
| | high | 0.92 | 11,022 | 24% | 0.89 | 10,671 | 23% | 0.93 | 15,612 | 37% | 0.88 | 13,323 | 28% | 0.91 | 10,375 | 22% | 0.87 | 10,077 | 21% |
| | **Ø** | **0.92** | **7,033** | **24%** | **0.87** | **6,832** | **23%** | **0.93** | **10,030** | **39%** | **0.70** | **13,368** | **32%** | **0.91** | **6,610** | **22%** | **0.86** | **6,446** | **21%** |

| | | production structure and demand scenarios | | | | | | | | | | | | | | | | |
|---|---|---|---|---|---|---|---|---|---|---|---|---|---|---|---|---|---|---|
| | | f m c | | | f m s | | | f l c | | | f l s | | | j m c | | | j m s | | |
| $c^b$ backorder cost rate | | R | β | % | R | β | % | R | β | % | R | β | % | R | β | % | R | β | % |
| | low | 0.93 | 5,756 | 20% | 0.91 | 6,011 | 21% | 0.94 | 7,863 | 31% | 0.89 | 8,046 | 29% | 0.92 | 5,429 | 19% | 0.89 | 5,726 | 19% |
| | med | 0.92 | 6,730 | 23% | 0.88 | 6,571 | 22% | 0.93 | 9,436 | 37% | 0.91 | 8,928 | 32% | 0.91 | 6,274 | 21% | 0.87 | 6,232 | 21% |
| | high | 0.91 | 8,613 | 29% | 0.83 | 7,915 | 27% | 0.92 | 12,791 | 50% | 0.90 | 10,542 | 35% | 0.88 | 8,128 | 27% | 0.81 | 7,379 | 24% |
| | **Ø** | **0.92** | **7,033** | **24%** | **0.87** | **6,832** | **23%** | **0.93** | **10,030** | **39%** | **0.70** | **13,368** | **32%** | **0.91** | **6,610** | **22%** | **0.86** | **6,446** | **21%** |

\* Sig. for R is ≤ 0.000

| | | | |
|---|---|---|---|
| f ... flow shop | j ... job shop | Ø ... average value | $c^i/c^e$ ... capacity cost rates internal / external |
| m ... many products | l ... low number of products | ρ ... shop load | $c^b$ ... backorder cost rate |
| c ... constant demand pattern | s ... seasonal demand pattern | | |

## 5 Conclusions

In this paper the effect of unbiased forecast errors on a hierarchical production planning system (MRP II) is investigated. The APP optimization problem is extended by including a planned utilization factor to mitigate the effect of assuming a deterministic setting in the APP but facing stochastic customer behaviour. A high cost reduction potential is identified when consciously applying this planned utilization factor as decision parameter. It is shown that for different production system structures, a set of different demand patterns and cost rates the planned utilization factor, and therefore the excess capacity, increases with respect to forecast error. The overall costs are as well found to significantly increase with forecast error (consistent over all parameter settings). For the basic scenario we find that an increase of the forecast error parameter α from 0.2 to 0.4 leads to a cost increase of approximately 10%. The results show that considerable cost reductions can be gained through improved collaboration of supply chain members on forecast accuracy. These findings on the effects of forecast error are in line with earlier studies. However, the contribution of this study lies in the focus on the long term forecasts provided and not on required lead times or order amounts stated by the customers. Note that the results are limited to situations with unbiased forecast errors. An additional contribution of this paper is the rigid modelling of the MRP II concept with the MPS and MRP step included which is shown to already mitigate some of the negative effects of short term information uncertainty. A further limitation of the results is that simulation is needed to identify the optimal planned utilization value. Nevertheless, a practical approach to identify a good planned utilization value can be suggested based on the observed convex overall cost function, if simulation is not feasible. Starting with a moderate planned utilization factor, this value can be increased until backorder costs start to become significant and overall costs increase. Then the value with the lowest overall cost or with a still acceptable service level can be chosen.

Further research could include the production system uncertainties concerning processing times, setup times and machine failures. Additionally, further investigation is required to study other demand functions, to investigate mis-specified forecasts (biased forecasts), and to analyse the impact of demand autocorrelation.

**Acknowledgments**

We are grateful to two anonymous referees for their comments that helped us notably to clarify the key concepts of this paper. Additionally, financial support from the Austrian Research Promotion Agency (FFG) under grant #826789 is gratefully acknowledged.

This is an Accepted Manuscript of an article published by Taylor & Francis in International Journal of Production Research on 25 Mar 2016, available at https://www.tandfonline.com/doi/abs/10.1080/00207543.2016.1162918.

**References**

Altendorfer, Klaus, and Stefan Minner. 2011. "Simultaneous optimization of capacity and planned lead time in a two-stage production system with different customer due dates." *European Journal of Operational Research* 213 (1): 134–46. doi:10.1016/j.ejor.2011.03.006.

Andersson, Håkan, Sven Axsäter, and Henrik Jönsson. 1981. "Hierarchical material requirements planning." *International Journal of Production Research* 19 (1): 45. doi:10.1080/00207548108956628.

Axsäter, S. 2000. *Inventory Control.* International Series in Operations Research & Management Science. Boston: Kluwer Academic Publishers.

Chern, C.-C., and J.-S. Hsieh. 2007. "A heuristic algorithm for master planning that satisfies multiple objectives." *Computers & Operations Research* 34 (11): 3491–3513. doi:10.1016/j.cor.2006.02.022.

Cui, Ruomeng, Gad Allon, Achal Bassamboo, and Jan A. van Mieghem. 2015. "Information Sharing in Supply Chains: An Empirical and Theoretical Valuation." *Management Science* 61 (11): 2803–24. doi:10.1287/mnsc.2014.2132.

Dolgui, Alexandre, and Caroline Prodhon. 2007. "Supply planning under uncertainties in MRP environments: A state of the art." *Annual Reviews in Control* 31 (2): 269–79. doi:10.1016/j.arcontrol.2007.02.007.

El-Bouri, Ahmed. 2012. "A cooperative dispatching approach for minimizing mean tardiness in a dynamic flowshop." *Computers & Operations Research* 39 (7): 1305–14. doi:10.1016/j.cor.2011.07.004.

Enns, S. T. 2001. "MRP performance effects due to lot size and planned lead time settings." *International Journal of Production Research* 39 (3): 461–80.

Enns, S. T. 2002. "MRP performance effects due to forecast bias and demand uncertainty." *European Journal of Operational Research* 138 (1): 87–102.

Felberbauer, Thomas, and Klaus Altendorfer. 2014. "Comparing the performance of two different customer order behaviors within the hierarchical production planning." In *Proceedings of the Winter Simulation Conference 2014*, 2227–38. Savannah GA.

Felberbauer, Thomas, Klaus Altendorfer, Gruber Daniel, and Alexander Hübl. 2013. "Application of a Generic Simulation Model to Optimize Production and Workforce Planning at an Automotive Supplier." In *Proceedings of the Winter Simulation Conference 2013*, 2689–97. Washington DC.

Felberbauer, Thomas, Klaus Altendorfer, and Alexander Hübl. 2012. "Using a scalable simulation model to evaluate the performance of production system segmentation in a combined MRP and kanban system." In *Proceedings of the Winter Simulation Conference 2012*, 1–12. Berlin.

Fildes, R., and B. Kingsman. 2010. "Incorporating demand uncertainty and forecast error in supply chain planning models." *J Oper Res Soc* 62 (3): 483–500. doi:10.1057/jors.2010.40.

Fleischmann, B., and H. Meyr. 2003. "Planning Hierarchy, Modeling and Advanced Planning Systems." In *Supply Chain Management: Design, Coordination and Operation*. Vol. 11, edited by S. C. Graves and A. G. de Kok, 455–523. Handbooks in Operations Research and Management Science: Elsevier.

Framinan, J. M., R. Ruiz-Usano, and R. Leisten. 2000. "Input control and dispatching rules in a dynamic CONWIP flow-shop." *International Journal of Production Research* 38 (18): 4589–98. doi:10.1080/00207540050205523.

Gansterer, Margaretha. 2015. "Aggregate planning and forecasting in make-to-order production systems." *International Journal of Production Economics.* doi:10.1016/j.ijpe.2015.06.001.

Hax, A. C., and H. C. Meal. 1975. "Hierarchical Integration of Production Planning and Scheduling." In *Logistics: TMS Studies in the Management Sciences*, edited by M. A. Geister, 53–69. North-Holland.

Ho, Chrwan-Jyh, and Tim C. Ireland. 1993. "A diagnostic analysis of the impact of forecast errors on production planning via MRP system nervousness." *Production Planning & Control* 4 (4): 311–22. doi:10.1080/09537289308919453.

Ho, Chrwan-Jyh, and Tim C. Ireland. 1998. "Correlating MRP system nervousness with forecast errors." *International Journal of Production Research* 36 (8): 2285–99. doi:10.1080/002075498192904.

Ho, Chrwan-Jyh, and Tim C. Ireland. 2012. "Mitigating forecast errors by lot-sizing rules in ERP-controlled manufacturing systems." *International Journal of Production Research* 50 (11): 3080–94. doi:10.1080/00207543.2011.592156.

Hopp, W. J., and M. L. Spearman. 2008. *Factory Physics:*. 3rd ed. Boston: Mc Graw Hill / Irwin.

Hübl, Alexander, Klaus Altendorfer, Herbert Jodlbauer, Margaretha Gansterer, and Richard F. Hartl. 2011. "Flexible model for analyzing production systems with discrete event simulation." In *Proceedings of the Winter Simulation Conference*, 1559–70. WSC '11: Winter Simulation Conference. http://dl.acm.org/citation.cfm?id=2431518.2431701.

Jeunet, Jully. 2006. "Demand forecast accuracy and performance of inventory policies under multi-level rolling schedule environments." *International Journal of Production Economics* 103 (1): 401–19. doi:10.1016/j.ijpe.2005.10.003.

Jodlbauer, Herbert, and Klaus Altendorfer. 2010. "Trade-off between capacity invested and inventory needed." *European Journal of Operational Research* 203 (1): 118–33. doi:10.1016/j.ejor.2009.07.011.

Jonsson, P., and S. Mattsson. 2006. "A longitudinal study of material planning applications in manufacturing companies." *International Journal of Operations & Production Management* 26 (9-10): 971–95. doi:10.1108/01443570610682599.

Kok, T. G. de, and J. C. Fransoo. 2003. "Planning Supply Chain Operations: Definition and Comparison of Planning Concepts." In *Supply Chain Management: Design, Coordination and Operation*. Vol. 11, edited by S.C. Graves and A.G. de Kok, 597–675. Handbooks in Operations Research and Management Science: Elsevier.

Kung, Ling-Chieh, and Ching-Chin Chern. 2009. "Heuristic factory planning algorithm for advanced planning and scheduling." *Computers & Operations Research* 36 (9): 2513–30. doi:10.1016/j.cor.2008.09.013.

Lee, T. S., and Everett E. Adam. 1986. "Forecasting Error Evaluation in Material Requirements Planning (MRP) Production-Inventory Systems." *Management Science* 32 (9): 1186–1205. doi:10.1287/mnsc.32.9.1186.

Louly, Mohamed-Aly, and Alexandre Dolgui. 2013. "Optimal MRP parameters for a single item inventory with random replenishment lead time, POQ policy and service level constraint." *International Journal of Production Economics* 143 (1): 35–40. doi:10.1016/j.ijpe.2011.02.009.

Meal, Harlan C. 1984. "Putting production decisions where they belong." *Harvard Business Review* 62 (2): 102–11.

Mieghem, J. A. Van. 1998. "Investment Strategies for Flexible Resources." *Management Science* 44 (8): 1071–78. doi:10.1287/mnsc.44.8.1071.

Missbauer, H., and R. Uzsoy. 2011. "Optimization Models of Production Planning Problems." In *Planning production and inventories in the extended enterprise: A state of the art handbook, volume 1*, edited by Karl G. Kempf, Pinar Keskinocak, and Reha Uzsoy, 437–507. New York: Springer.

Orlicky, J. 1975. *Material requirements planning;: The new way of life in production and inventory management.* New York: McGraw-Hill.

Panwalkar, S. S., and W. Iskander. 1977. "A Survey of Scheduling Rules." *Operations Research* 25 (1): 45–61.

Rajagopalan, S., and J. M. Swaminathan. 2001. "A Coordinated Production Planning Model with Capacity Expansion and Inventory Management." *Management Science* 47 (11): 1562–80. doi:10.1287/mnsc.47.11.1562.10254.

Rossi, Roberto, Onur A. Kilic, and S. A. Tarim. 2015. "Piecewise linear approximations for the static–dynamic uncertainty strategy in stochastic lot-sizing." *Omega* 50: 126–40.

Schneeweiss, Christoph. 2003. *Distributed decision making.* 2nd ed. Berlin [u.a.]: Springer.

Shafer, Scott M., and Timothy L. Smunt. 2004. "Empirical simulation studies in operations management: context, trends, and research opportunities." *Journal of Operations Management* 22 (4): 345–54. doi:10.1016/j.jom.2004.05.002.

Trapero, Juan R., N. Kourentzes, and R. Fildes. 2012. "Impact of information exchange on supplier forecasting performance." *Omega* 40 (6): 738–47. doi:10.1016/j.omega.2011.08.009.

Tsai, Shing C., and Chung H. Liu. 2015. "A simulation-based decision support system for a multi-echelon inventory problem with service level constraints." *Computers & Operations Research* 53: 118–27.

Wijngaard, J., and F. Karaesmen. 2007. "Advance demand information and a restricted production capacity: on the optimality of order base-stock policies." *OR Spectrum* 29 (4): 643–60. doi:10.1007/s00291-006-0076-x.

Yeung, J. H., W.C.K Wong, and L. Ma. 1998. "Parameters affecting the effectiveness of MRP systems: A review." *International Journal of Production Research* 36 (2): 313–32. doi:10.1080/002075498193750.